\documentclass[preprint,12pt]{elsarticle}
\usepackage{amsmath,amsfonts,graphics,epsfig,calc,color,array,url}

\newtheorem{theorem}{Theorem}

\newenvironment{proof}
{\noindent\textbf{Proof.}~\itshape}
{\par}

\usepackage{amssymb}

\makeatletter
\def\ps@pprintTitle{%
  \let\@oddhead\@empty
  \let\@evenhead\@empty
  \let\@oddfoot\@empty
  \let\@evenfoot\@oddfoot
}
\makeatother

\begin{document}

\begin{frontmatter}



\title{Approximating the inverse of a balanced symmetric matrix with positive elements}


\author{Ting Yan$^\dag$ and Xu Jinfeng$^\ddag$}

\address{\small $^\dag$Department of Statistics,
\small Central China Normal University, Wuhan, 430079, \\P. R. China\\
\small $^\ddag$ Department of Statistics and Applied Probability,
\small
National University of Singapore, 6 Science Drive 2, Singapore 117546, Singapore}

\begin{abstract}
For an $n\times n$ balanced symmetric matrix $T=(t_{i,j})$
with positive elements satisfying
$t_{i,i}= \sum_{j\neq i} t_{i,j}$ and certain
bounding conditions, we propose to use the matrix $S=(s_{i,j})$ to
approximate its inverse, where
$s_{i,j}=\delta_{i,j}/t_{i,i}-1/t_{..}$, $\delta_{i,j}$ is the
Kronecker delta function, and $t_{..}=\sum_{i,j=1
}^{n}(1-\delta_{i,j}) t_{i,j}$.  An explicit bound on the
approximation error is obtained, showing that the inverse is well
approximated to order $1/(n-1)^2$ uniformly.
\end{abstract}


\end{frontmatter}

\section{Introduction}
When solving a large system of linear equations, an accurate
approximation of the inverse of the coefficient matrix is crucially
important in establishing fast convergence rates for iterative
algorithms. For extensive reviews, see, for example,
\cite{Axelsson85, Benzi, Bruaset95, Zhang09}. In this paper, we
consider the approximation of the inverse of an $n\times n$ balanced
symmetric matrix $T=(t_{i,j})$ with positive
elements, i.e.,
\begin{equation}\label{eq1}
t_{i,j}=t_{j,i}>0 \mbox{~~and~~} t_{i,i}= \sum\limits_{j=1,j\neq i}^{n}t_{i,j},~~~ i=1,\cdots, n.
\end{equation}
The matrix $T$ is a special case of the diagonally dominant nonnegative matrix that has received wide
attention \cite{Dahl2000, Kaykobad, Maggs2002}. It is easy to show
that $T$ must be positive definite. The inverse of a general
nonnegative matrix has been extensively studied by \cite{Berman,
Loewya78, Eglestona2004, Markham, Martinez94}.

We propose to
approximate the inverse of $T$, $T^{-1}$, by the matrix
$S=(s_{i,j})$, where
\begin{equation*}
s_{i,j}=\frac{\delta_{i,j}}{t_{i,i}}-\frac{1}{t_{..}},
\end{equation*}
where $t_{..}=\sum_{i,j=1 }^{n}(1-\delta_{i,j}) t_{i,j}$. An
explicit upper bound on the approximation error is given in the
following section, which is crucially useful in establishing the
asymptotical normality of an estimated vector in the $\beta$-model
for undirected random graphs with a diverging number of nodes
\cite{Yan}.

\section{An explicit bound on the approximation error}
Let $m:=\min\limits_{1\le i<j\le n}
t_{i,j}$ and
$M:= \max\limits_{1\le i<j\le n} t_{i,j}$,
and for a matrix
$A=(a_{i,j})$, define $||A||:=\max_{i,j} |a_{i,j}|$. We have the
following theorem.
\begin{theorem}\label{dense-theorem}
$$||T^{-1}-S||\le \frac{C(m, M)}{(n-1)^2} + \frac{1}{2m(n-1)^2},$$
\end{theorem}
where
\begin{equation*}
C(m, M)=\frac{M}{m^2}\times [\frac{nM+(n-2)m }{ 2(n-2)m }].
\end{equation*}

\begin{proof}
Let $I_n$ be the $n\times n$ identity matrix. Define $F=T^{-1}-S$,
$V=(v_{ij})=I_n-TS$ and $W=(w_{ij})=SV$. We have the recursion
\begin{equation}\label{recursion}
F=T^{-1}-S=(T^{-1}-S)(I_n-TS)+S(I_n-TS)=FV+W.
\end{equation}
 Note that
\begin{eqnarray}
\nonumber
v_{i,j}&=&\delta_{i,j}-\sum_{k=1}^n t_{i,k}s_{k,j}\\
\nonumber
       &=&\delta_{i,j}-\sum_{k=1}^n t_{i,k}(\frac{\delta_{k,j}}{ t_{j,j} }-\frac{1}{t_{..}})\\
      \label{vij}
       &=&(\delta_{i,j}-1)\frac{t_{i,j}}{t_{j,j}}+\frac{2t_{i,i}}{t_{..}},
\end{eqnarray}
and
\begin{eqnarray}
\nonumber
w_{i,j}&=&\sum_{k=1}^n s_{i,k}v_{k,j}=\sum_{k=1}^n( \frac{\delta_{i,k}}{t_{i,i}}-\frac{1}{t_{..}} )
            [ (\delta_{k,j}-1)\frac{t_{k,j}}{t_{j,j}}+\frac{2t_{k,k}}{t_{..}}]\\
            \nonumber
      &=& \sum_{k=1}^n \frac{\delta_{i,k}}{t_{i,i}}  [(\delta_{k,j}-1)\frac{t_{k,j}}{t_{j,j}}+
          \frac{2t_{k,k}}{t_{..}}]-\frac{1}{t_{..}}\sum_{k=1}^n [(\delta_{k,j}-1)\frac{t_{k,j}}{t_{j,j}}
            +\frac{2t_{k,k}}{t_{..}}] \\
            \nonumber
      &=&[\frac{(\delta_{i,j}-1)}{t_{i,i}}( \frac{t_{i,j}}{t_{j,j}})+\frac{2t_{i,i}}{t_{i,i}t_{..}}]
         -\frac{1}{t_{..}}(\frac{-t_{j,j}}{t_{j,j}}+2 )\\
         \label{wij}
      &=&\frac{(\delta_{i,j}-1)t_{i,j}}{t_{i,i}t_{j,j}}+\frac{1}{t_{..}}.
\end{eqnarray}
Furthermore, when $i\neq j$,
\begin{eqnarray*}
&&0<\frac{1}{t_{..}}\le \frac{1}{mn(n-1)},\\
&&0<\frac{t_{i,j}}{t_{i,i}t_{j,j}}\le \frac{M}{m^2(n-1)^2},
\end{eqnarray*}
and it is easy to show, when $i,j,k$ are different from each other,
\begin{eqnarray*}
|w_{i,i}|&\le & \frac{1}{mn(n-1)},\\
|w_{i,j}|&\le & \frac{1}{mn(n-1)}, \\
|w_{i,j}-w_{i,k}| & \le & \frac{M}{m^2(n-1)^2}, \\
|w_{i,i}-w_{i,k}| & \le & \frac{M}{m^2(n-1)^2}.
\end{eqnarray*}
It follows that
\begin{equation}\label{wijin}
\max(|w_{i,j}|, |w_{i,j}-w_{i,k}|)\le \frac{M}{m^2(n-1)^2} ~~~~~ \mbox{for all $i,j,k$}.
\end{equation}
Next we use the recursion \eqref{recursion} to obtain a bound of the
approximate error $||F||$.  Let $a=\frac{M}{m^2(n-1)^2}$. By
\eqref{recursion} and \eqref{vij}, for any $i$, we have
\begin{equation}\label{fij}
f_{i,j}=\sum_{k=1}^n f_{i,k}[(\delta_{k,j}-1)\frac{t_{k,j}}{t_{j,j}}+\frac{2t_{k,k}}{t_{..}}] +w_{i,j}, ~~~~~j=1,\cdots,n.
\end{equation}
Thus, to prove Theorem 1, it is sufficient to show that
$|f_{i,j}|\le C(M,m)/(n-1)^2$ for any $i,j$. Fixing any
$i$, let $f_{i,\alpha}=\max\limits_{1\le k\le n}f_{i,k}$ and
$f_{i,\beta}=\min\limits_{1\le k\le n} f_{i,k}$.

First, we will show that $f_{i,\beta}\le 1/t_{..}\le 1/(m(n-1)^2)$.
A direct calculation gives that
\begin{eqnarray}  \nonumber
\sum_{k=1}^n f_{i,k}t_{k,i}& = & \sum_{k=1}^n (T_{i,k}^{-1}-(\frac{\delta_{i,k}}{t_{i,i}}-\frac{1}{t_{..}}) )t_{k,i}\\
\label{hold}
&=& 1 -  (1-\sum_{k=1}^n \frac{t_{k,i}}{t_{..}})=\sum_{k=1}^n \frac{t_{k,i}}{t_{..}}.
\end{eqnarray}
Thus, $f_{i,\beta} \sum_{k=1}^n t_{k,i} \le \sum_{k=1}^n f_{i,k}t_{k,i} = \sum_{k=1}^n \frac{t_{k,i}}{t_{..}}$. It
follows that $f_{i,\beta}\le 1/t_{..}$ and, similarly, $f_{i,\alpha}\ge 1/t_{..}$.

Note that $
f_{i,\beta}=-\sum_{k=1}^n f_{i,\beta}(\delta_{k,\alpha}-1)\frac{t_{k,\alpha}}{t_{\alpha,\alpha}}$. Thus,
\begin{equation}\label{falpha}
f_{i,\alpha}+f_{i,\beta}=\sum_{k=1}^n
(f_{i,k}-f_{i,\beta})(\delta_{k,\alpha}-1)\frac{t_{k,\alpha}}{t_{\alpha,\alpha}}
+\sum_{k=1}^n f_{i,k}(\frac{2t_{k,k}}{t_{..}})+w_{i,\alpha}.
\end{equation}
Similarly, we have that
\begin{equation}\label{fbeta}
f_{i,\beta}+f_{i,\beta}=\sum_{k=1}^n
(f_{i,k}-f_{i,\beta})(\delta_{k,\beta}-1)\frac{t_{k,\beta}}{t_{\beta,\beta}}
+\sum_{k=1}^n f_{i,k}(\frac{2t_{k,k}}{t_{..}})+w_{i,\beta}.
\end{equation}
Combining the above two equations, it yields
\begin{equation}\label{falphabeta}
 f_{i,\alpha}-f_{i,\beta}
= \sum_{k=1}^n (f_{i,k}-f_{i,\beta})[(\delta_{k,\alpha}-1)\frac{t_{k,\alpha}}{t_{\alpha,\alpha}}
-(\delta_{k,\beta}-1)\frac{t_{k,\beta}}{t_{\beta,\beta}}]+w_{i,\alpha}-w_{i,\beta} .
\end{equation}
Let $\Omega=\{k:(1-\delta_{k,\beta})t_{k,\beta}/t_{\beta,\beta}
\ge (1-\delta_{k,\alpha})t_{k,\alpha}/t_{\alpha,\alpha} \}$ and let $|\Omega|=\lambda$.
Note that $1\le \lambda \le n-1$. Then,
\begin{eqnarray}\nonumber
&&\sum_{k=1}^n(f_{i,k}-f_{i,\beta})[(\delta_{k,\alpha}-1)\frac{t_{k,\alpha}}{t_{\alpha,\alpha}}
-(\delta_{k,\beta}-1)\frac{t_{k,\beta}}{t_{\beta,\beta}}] \\
\nonumber &\le &\sum_{k\in \Omega}(f_{i,k}-f_{i,\beta})[
(1-\delta_{k,\beta})\frac{t_{k,\beta}}{t_{\beta,\beta}}-(1-\delta_{k,\alpha})\frac{t_{k,\alpha}}{t_{\alpha,\alpha}}] \\
\nonumber
&\le &(f_{i,\alpha}-f_{i,\beta}) [\frac{\sum_{k\in
\Omega}t_{k,\beta}}{t_{\beta,\beta}}-\frac{\sum_{k\in
\Omega}(1-\delta_{k,\alpha})t_{k,\alpha}}{t_{\alpha,\alpha}}] \\
\nonumber
 & \le & (f_{i,\alpha}-f_{i,\beta}) [\frac{\lambda M
}{\lambda M+(n-1-\lambda)m}-\frac{(\lambda-1) m}{(\lambda-1)
m+(n-\lambda)M}].\\
\label{yyy}&&
\end{eqnarray}
Let
\begin{equation*}
f(\lambda)=\frac{\lambda M}{\lambda M+ (n-1-\lambda)m }-
\frac{(\lambda-1)m}{ (\lambda-1)m+ (n-\lambda)M}.
\end{equation*}
There are two cases to consider the maximum of
$f(\lambda)$ in the range of
$\lambda\in [1, n-1]$.\\
Case I: When $M=m$, it is easy to show $f(\lambda)=1/(n-1)$.\\
Case II: $M\neq m$.
Since
\renewcommand\arraystretch{1.5}
\begin{eqnarray*}
\begin{array}{lll}
f^\prime(\lambda) &=& \frac{ (n-1)Mm}{[\lambda M+
(n-1-\lambda)m]^2 } -
\frac{ (n-1)Mm}{ [(\lambda-1)m+(n-\lambda)M]^2} \\
&=& \frac{ (n-1)Mm [(n-2\lambda)(M-m)] [\lambda
M+ (n-1-\lambda)m+(\lambda-1)m+ (n-\lambda)M]}{[\lambda M+
(n-1-\lambda)m]^2[(\lambda-1)m+ (n-\lambda)M]^2 }
\end{array}
\end{eqnarray*}
and
\begin{equation*}
\begin{array}{c}
f^{\prime\prime}(\lambda)=-2(M-m)Mm(n-1)\left(\frac{ 1}{[\lambda M+
(n-1-\lambda)m]^3 } + \frac{ 1}{ [(\lambda-1)m+
(n-\lambda)M]^3} \right),
\end{array}
\end{equation*}
$f(\lambda)$ takes its maximum at $\lambda=n/2$ when $1\le \lambda\le n-1$.  A direct calculation gives that
\begin{eqnarray}\label{flambdamax}
f(\frac{n}{2}) & = & \frac{ nM-(n-2)m}{ nM+(n-2)m}.
\end{eqnarray}
Combining \eqref{falphabeta}, \eqref{yyy} and \eqref{flambdamax}, it yields
\begin{eqnarray*}
f_{i,\alpha}-f_{i,\beta}
\le   \frac{nM-(n-2)m}{ nM+(n-2)m}\times (f_{i,\alpha}-f_{i,\beta}) + a,
\end{eqnarray*}
so that
\begin{equation*}
f_{i,\alpha}-f_{i,\beta}
 \le    [\frac{2(n-2)m }{ nM+(n-2)m}]^{-1}   \times a =C(M, m)/(n-1)^2.
\end{equation*}
Note that
\[
\max_{j=1,\cdots, n}|f_{i,j}| \le   f_{i,\alpha}-f_{i,\beta}+\frac{1}{t_{..}} \le
[\frac{nM+(n-2)m }{2(n-2)m }]\times \frac{M}{m^2(n-1)^2}+\frac{1}{2m(n-1)^2}.
\]
This completes the proof.
\end{proof}

\section{Discussion}
In many applications, it is important to closely approximate the
inverse of a matrix when its explicit form is unavailable. For
example, when an algorithm involves solving a matrix in each
iteration, its convergence rate is often related to the approximate
inverse it uses. On the other hand, in some statistical
applications, an accurate approximation of the inverse of the Fisher
information matrix is critical in establishing the theoretical
properties of the maximum likelihood estimates. For instance, Simons
and Yao [15] obtained a good approximation of the inverse of a
symmetric positive definite matrix with negative off-diagonal
elements. This result is crucial in establishing their most
surprising result that the maximum likelihood estimates of the merit
parameters in the Bradley-Terry model for paired comparisons retain
good asymptotic properties even when the number of subjects goes to
infinity. Similarly, our results can be readily used to prove the
asymptotic normality of the maximum likelihood estimate in the
$\beta$-model with a diverging dimension \cite{Yan} since the Fisher
information matrix of the $\beta$ model is a diagonally dominant nonnegative
matrix.

The matrix $S$ which we use to approximate the inverse of the matrix
$T$ takes the form of $I+H_c$, where each element of $H_c$ is $c$.
If $c>0$, then $S$ is a class of preconditioners for $M$-matrices
\cite{Zhang09}. In our  situation, $c<0$ since $S$ is a matrix with
non-negative elements. The bound on the approximation error  in
Theorem 1  depends on $m$, $M$ and $n$. When $m$ and $M$ are bounded
by a constant, all the elements of $T^{-1}-S$ are of order
$O(1/(n-1)^2)$ as $n\to\infty$, uniformly.

We illustrate by an example that the bound on the
approximation error in Theorem 2.1 is optimal in the sense that any bound in the
form of $C(m,M)/f(n)$ requires $f(n)=O((n-1)^2)$ as $n\to\infty$.
Assume that the matrix $T$ consists of the elements:
$t_{i,i}=(n-1)M, i=1,\cdots,n-1; t_{n,n}=(n-1)m$ and
$t_{i,j}=m, i,j=1,\cdots, n; i\neq j$, which satisfies \eqref{eq1}. By the Sherman-Morrison formula, we
have
\begin{eqnarray*}
\scriptsize
\begin{array}{lll}
(T^{-1})_{i,j}&=&\frac{\delta_{i,j} }{(n-1)M-m}-\frac{m}{[(n-1)M-m]^2}\times \left(1+\frac{(n-1)m}{(n-1)M-m}+
\frac{1}{(n-2)}\right)^{-1}, i,j=1,\cdots,n-1\\
(T^{-1})_{n,j}&=&\frac{\delta_{n,j} }{(n-2)m}-\frac{1}{(n-2)[(n-1)M-m]}\times \left(1+\frac{(n-1)m}{(n-1)M-m}+
\frac{1}{(n-2)}\right)^{-1},~~j=1,\cdots,n-1\\
(T^{-1})_{n,n}&=&\frac{1}{(n-2)m}-\frac{1}{(n-2)^2m}\times \left(1+\frac{(n-1)m}{(n-1)M-m}+
\frac{1}{(n-2)}\right)^{-1}.
\end{array}
\end{eqnarray*}
In this case, the elements of $S$ are
\begin{eqnarray*}
S_{i,j} & = &\frac{\delta_{i,j}}{(n-1)M}-\frac{1}{n(n-1)m}, ~~~i,j=1,\cdots, n-1;i\neq j,\\
S_{n,j} & = &\frac{\delta_{n,j}}{(n-1)m}-\frac{1}{n(n-1)m}, ~~~j=1,\cdots, n.
\end{eqnarray*}
It is easy to show that the bound of $||T^{-1}-S||$ is $
\frac{1}{(n-1)^2m} +o(\frac{1}{(n-1)^2})$. This suggests that the rate $1/(n-1)^2$ is
optimal. On the other hand, there is a gap between $1/m$ and
$C(m,M)=\frac{(M+m)M}{2m^3}+\frac{1}{m}+o(1)$ which implies that there might be space for
improvement.

Finally, we discuss some extension to diagonally dominant case, where we still use $S$ to approximate
the inverse of $T$. In this case, some $t_{i,i}$
may be greater than the corresponding row sum without the diagonal element. Denote $\Delta_i:= t_{i,i}-\sum_{j\neq i} t_{i,j} $
and redefine  $m:=\min\limits_{1\le i<j\le n}t_{i,j}$ and
$M:=\max\{ \max\limits_{1\le i<j\le n} t_{i,j}, \max\limits_{1\le i\le n} \Delta_i\}$. With a similar argument,
we can prove the new upper bound of the approximation errors is
\begin{scriptsize}
\[
\frac{1}{(n-1)^2}\times [\frac{2(n-2)m }{ nM+(n-2)m}-\frac{M}{m(n-1)}-\frac{ (n-2)Mm }{ [(n-2)m+M][ (n-2)m + 2M ] }]^{-1}
 \times (\frac{M}{m^2}+ \frac{4M}{m^2n})+\frac{1}{mn(n-1)},
\]
\end{scriptsize}
{\normalsize for large $n$ when $M/m=o(n)$.}







\bibliographystyle{model1a-num-names}



\end{document}